%% file: Poisson-like-va.tex
\documentclass[paper=A4,fontsize=12pt]{scrartcl}

\newcommand{\kmcomment}[1]{}

\input{km_cmd.tex}

\input{km_cmd2.tex}

\input{km-env-e.tex}

\usepackage[backend=biber,style=numeric]{biblatex}
\bibliography{km_refs.bib}

\usepackage[leqno]{amsmath}
\usepackage{amssymb,amsfonts}

\usepackage{comment,tabularray}
\UseTblrLibrary{amsmath,booktabs,counter,diagbox,siunitx,varwidth}

\usepackage{tikz}
\usetikzlibrary{positioning,calc,quotes,angles,cd}
\usetikzlibrary{backgrounds}
\usepackage{adjustbox}

\usepackage{newtxmath} \usepackage{newtxtext}
\usepackage{mathtools}
\usepackage{fancyvrb,shortvrb}
\usepackage{listings,jlisting,setspace,Xentry}

 \RequirePackage{listings} 
 \lstnewenvironment{kmverb}
    { \lstset{xleftmargin=10mm, breaklines=true}}
    { }

 \lstnewenvironment{kmverbb}[1]
    { \lstset{xleftmargin=#1, showspaces=false, breaklines=true}}
    { }
 \lstnewenvironment{kmverbN10}
    { \lstset{xleftmargin=-10mm, breaklines=true}}
    { } 
 \lstnewenvironment{kmverbZero}
    { \lstset{xleftmargin=0mm, breaklines=true}}
    { }

 \lstnewenvironment{kmverb*}
    { \lstset{xleftmargin=10mm, showspaces=true, breaklines=true}}
    { } 

\usepackage{wrapfig}
\usepackage{geometry}
\usepackage[leqno]{amsmath}
\usepackage{amssymb,amsfonts}
\usepackage{mathtools}
\usepackage{newtxmath}

\newcommand{\cb}[1]{\ANYb{c}{#1}}

\newcommand{\Cb}[1]{\ANYb{C}{#1}}

\DefineShortVerb{\!} 

\usepackage[all,warning]{onlyamsmath}
\numberwithin{equation}{section}

\usepackage{mathtools}

\usepackage{tikz}
\usetikzlibrary{positioning,calc, quotes, angles, decorations.pathmorphing}
\usetikzlibrary{backgrounds}

\makeatletter
\global\let\tikz@ensure@dollar@catcode=\relax
\makeatother

\usepackage{wrapfig}

\usepackage{amssymb,fge}

\title{Poisson-like cohomologies associated with some Lie superalgebras}

\author{
Kentaro Mikami
\thanks{Akita University 
(This work is supported by JSPS KAKENHI Grant Number 22K03306.)} 
 \and Tadayoshi Mizutani\thanks{Saitama University}
 \and Hajime Sato\thanks{Nagoya University}
 }

\allowdisplaybreaks
\AtBeginDocument{\parindent=0pt \setlength{\headheight}{16pt} 
} 
\begin{document}
\maketitle
\tableofcontents
\thispagestyle{plain}

\section{Introduction}
First, we recall briefly the Poisson cohomology of a Poisson manifold $M$.
A Poisson tensor \(\pi \) is a 2-vector field on \(M\) by which 
\(\Bkt{f}{g} = \langle \pi, df \we dg \rangle \) becomes the Poisson bracket
of functions \(f\) and \(g\) on $M$.  This is equivalent to 
\(\SbtS{\pi}{\pi} = 0\), where \(\SbtS{\cdot}{\cdot}\) is the
Schouten(-Nijenhuis) bracket on multi-vector fields. 
A quick introduction of the Schouten bracket is 
\begin{equation}
\SbtS { X_{1} \we  \cdots \we X_{r}}
{ Y_{1}  \we \cdots \we Y_{s}} = \sum \parity{i+j} \Sbt{X_{i}}
{Y_{j}} \we ( X_{1} \we \cdots \widehat{X_{i}} \cdots \we X_{r}) \we 
( Y_{1} \we \cdots \widehat{Y_{j}}  \cdots \we Y_{s})
\end{equation}
where \( \Sbt{X_{i}}{Y_{j}}\) is the Jacobi-Lie bracket of two vector
fields \(X_{i}\) and \(Y_{j}\).  In short,  
the Schouten bracket is a Lie superalgebra bracket on \(\sum_{j=1}^{\dim M} \Lambda^{j} \tbdl{M}\), where   
the grading of \(\Lambda^{j}\tbdl{M} = j-1\).

For a Poisson tensor \(\pi\), we define 
\( d_{\pi}(U) = \SbtS{\pi}{U}\) for \(\ds U \in \Lambda ^{j} \tbdl{M}\).
Then  \(\ds d_{\pi}\left( \Lambda ^{j} \tbdl{M}\right) \subset 
 \Lambda ^{j+1} \tbdl{M}\),  and 
\(\ds d_{\pi} \circ d_{\pi} = 0\) holds by super bracket property of the
Schouten bracket. Thus, \(d_{\pi}\) is a
coboundary operator of degree \(1\) and 
we obtain the so-called Poisson cohomology groups. 
We easily expect that an even \(p\)-multivector field 
 \(\Pi \in  \Lambda ^{p} \tbdl{M}\) with \(\SbtS{\Pi}{\Pi} = 0\) defines  
\( d_{\Pi}(U) = \SbtS{\Pi}{U}\), satisfies   
\(\ds d_{\Pi}\left( \Lambda ^{j} \tbdl{M}\right) \subset 
 \Lambda ^{j+p-1} \tbdl{M}\)  and 
\(\ds d_{\Pi} \circ d_{\Pi} = 0 \). 
In this note, 
we translate notion of Poisson cohomology groups of Poisson manifolds to 
those of corresponding concepts of Lie superalgebras.  
Each vector field \(X\) on \(M\) is a section of \(\tbdl{M}\) and 
exactly written as 
\(X\in \Gamma \tbdl{M}\), but here we denote it simply as \(X \in
\tbdl{M}\). This abbreviation also applies to all multi-vector fields and 
differential forms.  

\section{Poisson-like cohomologies of Lie superalgebra}

The abstract definition of a Lie superalgebra is the following:    
\begin{definition}[Lie superalgebra] \label{def:superalg}
Suppose a real vector space 
$\frakg $ is graded by integers  as 
\(\ds \frakg = \sum_{j\in \mZ} \frakg_{j} \)
and there exits a \(\mR\)-bilinear operation \(\Sbt{\cdot}{\cdot}\)  satisfying 
\begin{align}
& \Sbt{ \frakg_{i}}{ \frakg_{j}} \subset \frakg_{i+j} \label{cond:1} \\
& \Sbt{X}{Y} + (-1) ^{ x y} 
 \Sbt{Y}{X} = 0 \quad \text{ where }  X\in \frakg_{x} \text{ and } 
 Y\in \frakg_{y}  \\
& 
(-1)^{x z} \Sbt{ \Sbt{X}{Y}}{Z}  
+(-1)^{y x} \Sbt{ \Sbt{Y}{Z}}{X}  
+(-1)^{z y} \Sbt{ \Sbt{Z}{X}}{Y}  = 0\;.  
\label{super:Jacobi}
\end{align}
Then we call \(\frakg\) a \textit{\(\mZ\)-graded} (or \textit{pre}) Lie superalgebra.    
\end{definition}

\begin{remark} 
Super Jacobi identity \eqref{super:Jacobi} above is equivalent to one of
the following two.  
\begin{align}
\Sbt{ \Sbt{X}{Y}}{Z} &=  
\Sbt{X}{ \Sbt{Y}{Z} } + (-1)^{y z} \Sbt{ \Sbt{X}{Z}}{Y} \\  
\Sbt{X}{\Sbt{Y}{Z}} &=  
\Sbt{\Sbt{X}{Y}}{Z}  + (-1)^{x y} \Sbt{Y}{ \Sbt{X}{Z}} 
\end{align} 
\end{remark}

We translate and generalize the notion of Poisson cohomology groups above on
Lie superalgebras, and we call the generalized notion to 
\textit{Poisson-like} cohomology groups. 

\begin{prop}
\label{prop:Poisson-like-cohom}
Let \(\frakg = \sum_{i} \frakgN{i}\) be a Lie superalgebra. 
For odd \(p\), assume 
\( \pi \in \frakgN{p}\) 
satisfies \(\Sbt{\pi}{\pi} = 0\). 
Then the map \(d_{\pi}\) defined by  
\begin{align}
\label{def:d:pi}
& d_{\pi} (U) = \Sbt{\pi}{U} 
\qquad (U\in \frakg) 
\\\shortintertext{ satisfies }
& d_{\pi} (\frakgN{i}) \subset \frakgN{i+p} \quad \text{and}\quad  
  d_{\pi} \circ d_{\pi} = 0 \; . 
\end{align}
\end{prop}
\textbf{Proof:} For \( \pi  \in \frakgN{p} \), 
Jacobi identity shows 
\( \Sbt{\pi }{\Sbt{\pi }{U}} = \Sbt{\Sbt{\pi }{\pi }}{U} + \parity{p^2}\Sbt{\pi }{\Sbt{\pi }{U}}\).  
Since \(p\) is odd,  
\(2 \Sbt{\pi }{\Sbt{\pi }{U}} = \Sbt{\Sbt{\pi }{\pi }}{U}\) holds, and the assumption 
\(\Sbt{\pi }{\pi } =0\) implies 
\(2 \Sbt{\pi }{\Sbt{\pi }{U}} = 0 \). \kmqed 
\begin{remark}
We may call 
\(
d_{\pi}\) defined by \eqref{def:d:pi} the coboundary operator associated
with \(\pi\), or Poisson-like coboundary operator, simply.   
\end{remark}
\begin{remark}
Contrary to Proposition \ref{prop:Poisson-like-cohom} above, we start from \(P\in \frakgN{p} \) where
\(p\) is even. Then  
\( \Sbt{P} {P} =0 \) holds automatically because 
\(\Sbt{P}{P} + \parity{p^2} 
\Sbt{P} {P} =0 \)  
by symmetric property of super bracket and \(p\) being even.  
Since
\[(d_{P} \circ d_{P})(U) = \Sbt{ P}{\Sbt{P}{U}} = \Sbt{ \Sbt{P}{P} }{U} + 
\parity{p^2}
\Sbt{P}{\Sbt{P}{U}} = 0 + 
\Sbt{P}{\Sbt{P}{U}} = 
(d_{P} \circ d_{P})(U) \;, \] 
there is no information about \( d_{P} \circ d_{P}
\) in general for each even grade \(P\).   
We refer to Example \ref{exam:BadNews} for a concrete example. 
\end{remark}

\begin{defn} [generalized Poisson cohomology] 
In a more general context than Proposition 
\ref{prop:Poisson-like-cohom}, 
let \(\Phi\) be a linear operator on 
a Lie superalgebra
\(\frakg = \sum_{i} \frakgN{i}\) into itself satisfying 
\begin{equation} \label{ddZero}
   \Phi (\frakgN{i}) \subset \frakgN{i+p} \quad \text{and}\quad  
   \Phi  \circ \Phi = 0 \; . 
\end{equation}
We call \(\Phi\) the coboundary operator of degree \(p\). 
Then \(\Phi\) defines 
\begin{equation} \label{degreeP}
\mH^{i,p} = 
\ker( \Phi : \frakgN{i} \to \frakgN{i+p})/ \Phi(\frakgN{i-p})\;. 
\end{equation}
We call \eqref{degreeP} Poisson-like cohomology groups defined by \(\Phi\). 
\end{defn}

If \( \frakgN{i}\) are finite dimensional, 
the dimension of \(\mH^{i,p} \) makes sense and we may call them Betti
numbers. By adding one more assumption, 
\( \dim \frakgN{i} = 0\) for \( |i|> b\), i.e., 
the length of Lie superalgebra is finite, then we get the following result.  

\newcommand{\absp}{|p|} 

\begin{thm}
Consider a 
Lie superalgebra 
\(\frakg = \sum_{ |i| \leqq b }  \frakgN{i}\)  and 
each \(\frakgN{i}\) is finite dimensional.  Let \(\Phi\) be a linear map
satisfying \eqref{ddZero} in the definition above, i.e., a coboundary
operator of degree \(p\). 
Then we have an integer 
\begin{equation} \label{altBetti}
\sum_{k} \parity{k} \sum_{k \absp \leqq i < (k+1) \absp} \dim
\left(
\ker( \Phi : \frakgN{i} \to \frakgN{i+p})/ \Phi(\frakgN{i-p})\right) \;. 
\end{equation}
The above \eqref{altBetti} 
is independent of the operator \(\Phi\) and  
agrees with the value 
\begin{equation} \label{MMS}
\sum_{k} \parity{k} \sum_{k \absp \leqq i < (k+1) \absp} \dim \frakgN{i} \;. 
\end{equation}
\end{thm}
\textbf{Proof:}
\begin{align*}
&\quad \text{ The  \eqref{altBetti} above } \\
& = \sum_{k} \parity{k} \sum_{k \absp \leqq i < (k+1) \absp}(\dim
\ker( \Phi : \frakgN{i} \to \frakgN{i+p}) - \dim(\Phi(\frakgN{i-p}))
\\
\shortintertext{
by the rank theorem of linear map, we see 
\( \dim \frakgN{i} - \dim
\ker( \Phi : \frakgN{i} \to \frakgN{i+p}) = \dim \Phi(\frakgN{i}) 
\) and we get 
}
&= 
\sum_{k} \parity{k} \sum_{k \absp \leqq i < (k+1) \absp}(
 \dim \frakgN{i} -  \dim \Phi(\frakgN{i}) 
- \dim \Phi(\frakgN{i-p}) 
\\
&= 
\sum_{k} \parity{k} \sum_{k \absp \leqq i < (k+1) \absp}
 \dim \frakgN{i} 
- \sum_{k} \parity{k} \sum_{k \absp \leqq i < (k+1) \absp}
   \dim \Phi(\frakgN{i}) 
- \sum_{k} \parity{k} \sum_{k \absp \leqq i < (k+1) \absp}
 \dim \Phi(\frakgN{i-p})
\\
\shortintertext{
we replace \(i-p\) by \(j\), we get
}
&= 
\sum_{k} \parity{k} \sum_{k \absp \leqq i < (k+1) \absp}
 \dim \frakgN{i} 
- \sum_{k} \parity{k} \sum_{k \absp \leqq i < (k+1) \absp }
   \dim \Phi(\frakgN{i}) 
- \sum_{k} \parity{k} \sum_{k \absp \leqq p+j < (k+1) \absp}
 \dim \Phi(\frakgN{j})
\\
\shortintertext{we may put \(p = \varepsilon \absp\) where \(\varepsilon = \pm 1\), and} 
&= 
\sum_{k} \parity{k} \sum_{k \absp \leqq i < (k+1) \absp}
 \dim \frakgN{i} 
- \sum_{k} \parity{k} \sum_{k \absp \leqq i < (k+1) \absp}
   \dim \Phi(\frakgN{i}) 
- \sum_{k} \parity{k} \sum_{(k-\varepsilon) \absp \leqq j < (k-\varepsilon+1)
\absp }
 \dim \Phi(\frakgN{j})
\\
&= 
\sum_{k} \parity{k} \sum_{k \absp \leqq i < (k+1) \absp}
 \dim \frakgN{i} 
- \sum_{k} \parity{k} \sum_{k \absp \leqq i < (k+1) \absp}
   \dim \Phi(\frakgN{i}) 
- 
\sum_{k} \parity{k+\varepsilon} \sum_{k \absp \leqq j < (k+1) \absp}
 \dim \Phi(\frakgN{j})
 \\
&= 
\sum_{k} \parity{k} \sum_{k \absp \leqq i < (k+1) \absp}
 \dim \frakgN{i} 
- \sum_{k} \parity{k} \sum_{k \absp \leqq i < (k+1) \absp}
   \dim \Phi(\frakgN{i}) 
+ \sum_{k} \parity{k} \sum_{k \absp \leqq j < (k+1) \absp}
 \dim \Phi(\frakgN{j})
 \\
&= 
\sum_{k} \parity{k} \sum_{k \absp \leqq i < (k+1) \absp}
 \dim \frakgN{i} \; . 
\end{align*} 

\begin{remark}
When \(p=1\), \eqref{MMS} is just the Euler number, and so \eqref{MMS} is a
generalization of the Euler number.  
\end{remark}

\subsection{Examples of Poisson-like cohomologies of  superalgebra by 
a Lie algebra}

\newcommand{\uType}{\text{Type}}
Although we have a typical superalgebra \( \sum_{p=1}^{\dim M} \Lambda^{p}
\tbdl{M} \) for a manifold \(M\), their chain spaces are infinite
dimensional,  and it is hard to manipulate  those (co)homology groups.
If \(M\) is a finite dimensional Lie group, then its invariant multivector
fields are finite dimensional, and provides a superalgebra, which is the
direct sum of \(\frakgN{i} =
\Lambda^{i+1}\frakm \) for \(i=0,\ldots, n - 1\), where \(M\) is a
$n$-dimensional Lie group whose Lie algebra is \(\frakm\). 

\begin{exam} \label{Eg:Type1}
We focus on the 4-dimensional Lie algebra spanned by 
\( \yb{i} \) for \(i=1..4\) satisfying 
(minimum) bracket relations 
\( 
\uType[1] =\{
            [\yb{2},\yb{4}] = \yb{1},  
            [\yb{3},\yb{4}] = \yb{2} 
\} \) 
which is listed  in \cite{MR3030954}, and  
take a 2-vector field \[\pi =
 c_{1} \yb{1} \we \yb{2}  
+c_{2} \yb{1} \we \yb{3}  
+c_{3} \yb{1} \we \yb{4}  
+c_{4} \yb{2} \we \yb{3}  
+c_{5} \yb{2} \we \yb{4}  
+c_{6} \yb{3} \we \yb{4}  
\] where \( c_{j}\) are real constants.  
Since 
\[ \SbtS{\pi}{\pi}/2 = 
 (-\cb{2} \cb{6} + \cb{4} \cb{5})  \mathit{\&^{\we}} (\yb{1}, \yb{2}, \yb{3})
 + (-\cb{3} \cb{6} + \cb{5}^2) \mathit{\&^{\we}} (\yb{1}, \yb{2}, \yb{4}) +
 \cb{5} \cb{6} \mathit{\&^{\we}} (\yb{1}, \yb{3}, \yb{4}) + \cb{6}^2
 \mathit{\&^{\we}} (\yb{2}, \yb{3}, \yb{4}) \; , 
\] 
Poisson condition \(\SbtS{\pi}{\pi} = 0\) implies  
\[
\left\{- c_{2} c_{6}+ c_{4} c_{5}=0, - c_{3} c_{6}+ c_{5}^{2}=0, 
 c_{5} c_{6}=0,  c_{6}^{2}=0 \right\} .
\] 
Just solving Poisson condition, we get  
\( c_{1} ,c_{2} ,c_{3} ,c_{4}\) are free, and \( c_{5} {=} 0,c_{6} {=}
0\) . Thus,   
 \[\pi = \yb{1} \we ( c_{1} \yb{2}  +c_{2} \yb{3}  +c_{3} \yb{4}  ) 
+c_{4} \yb{2} \we \yb{3} \; ,   
\text{ and } \;   
\pi \we \pi = 
2 c_{3} c_{4} {\&}^{\we} \! \left(y_{1},y_{2},y_{3},y_{4}\right)
\; . \] 
Since 
\begin{align*}
& d_{\pi} [ \yb{1},\ldots, \yb{4} ] = 
 [d_{\pi} (\yb{1}), \ldots, d_{\pi}(\yb{4})] = 
\left[0,0,-c_{3} y_{1}{\we }y_{2},c_{4} y_{1}{\we}y_{3}
+c_{2} y_{1}{\we}y_{2}\right]
\\
& 
d_{\pi} [ \yb{1}\we \yb{2}, \ldots, \yb{3} \we \yb{4} ] = 
\left[0,0,0,0,c_{4} {\&^{\we}} \! \left(y_{1},y_{2},y_{3}\right),
-c_{2} {\&^{\we}} \! \left(y_{1},y_{2},y_{3}\right)
-c_{3} {\&^{\we}} \! \left(y_{1},y_{2},y_{4}\right)
\right]
\\
&
d_{\pi} [ \yb{1}\we \yb{2} \we \yb{3}, \ldots, \yb{2} \we \yb{3} \we \yb{4} ] = 
\left[0,0,0,0\right]\;, 
\end{align*}
their matrix representations of \(d_{\pi}\) are
\[ A(1,2) = 
\left[\begin{array}{cccccc}
0 & 0 & 0 & 0 & 0 & 0 
\\
 0 & 0 & 0 & 0 & 0 & 0 
\\
 -c_{3} & 0 & 0 & 0 & 0 & 0 
\\
 c_{2} & c_{4} & 0 & 0 & 0 & 0 
\end{array}\right], 
\quad 
A(2,3) = 
\left[\begin{array}{cccc}
0 & 0 & 0 & 0 
\\
 0 & 0 & 0 & 0 
\\
 0 & 0 & 0 & 0 
\\
 0 & 0 & 0 & 0 
\\
 c_{4} & 0 & 0 & 0 
\\
 -c_{2} & -c_{3} & 0 & 0 
\end{array}\right], 
\quad A(3,4) = 
\left[\begin{array}{c}
0 
\\
 0 
\\
 0 
\\
 0 
\end{array}\right] ,  
\]
so their ranks are 
\( \rank A(1,2) = \rank A(2,3) = \rank \begin{bmatrix} \cb{2} & \cb{3} \\
\cb{4} & 0 \end{bmatrix}\),  and \(\rank A(3,4) = 0\).

If \( c_{3} c_{4} \ne 0\),  i.e., if   
\(\pi\) is a ``symplectic structure'', then  
\( \rank A(1,2) = \rank A(2,3) =2, \rank A(3,4) = 0\) and the Betti numbers
are (2,2,2,1).  

If \( c_{3}= c_{4} = 0\) and \( c_{2} = 0\), then \(\pi = c_{1} \yb{1} \we
\yb{2}\) and  
\( \rank A(1,2) = \rank A(2,3) = \rank A(3,4) = 0\) and the Betti numbers
are (4,6,4,1).  

If \( c_{3}= c_{4} = 0\) and \( c_{2} \ne 0\), 
or 
if \( c_{3}  c_{4} = 0\) and \( c_{3}^2 + c_{4}^2 > 0\), then 
then 
\( \rank A(1,2) = \rank A(2,3) = 1\) and \( \rank A(3,4) = 0\) and the Betti numbers
are (3,4,3,1).  
\end{exam}
\begin{remark}
The Euler number is always \(1\) because 
\( \frakgN{0} = \mR^{4}\) and 
\( \frakgN{j} = \Lambda^{j+1}\frakgN{0}\) for \(j=0,\ldots, 3 \), but  
\( \Lambda^{0}\frakgN{0} = \mR \) is not counted. 
\end{remark}
\begin{remark}
\cite{MR3030954} gives us 12 different types of 4-dimensional Lie algebra.
Type[1] has been introduced above.  For instance, 
the second one is given by 
\( 
\uType[2] =\{
            [\yb{1},\yb{4}] =  a   \yb{1},  
            [\yb{2},\yb{4}] = \yb{2}, 
            [\yb{3},\yb{4}] = \yb{2} + \yb{3} \} 
\). It has a parameter \(a\), the Poisson condition has five cases and so the
story becomes a little complicated, but still the same discussion works well
as in Type[1].     
\end{remark}

\section{Poisson-like cohomology of superalgebra of differential forms}
Let \(M\) be a differential manifold.  
The direct sum \(\sum_{j=0}^{\dim M} \Lambda^{j}\cbdl{M}\) 
with the bracket 
\begin{equation}
\Sbt{\alpha}{\beta} = \parity{a} d( \alpha \we \beta )\quad 
\text{where }\quad \alpha \in \Lambda^{a} \cbdl{M} 
\end{equation} 
forms a Lie superalgebra with the grading of \(a\)-form to be \(-(a+1)\) 
(cf. \cite{Mik:Miz:superForms}). However, another 
grading of \(a\)-form by \(a+1\) makes sense.  It is clear the identity map gives
Lie superalgebra isomorphism between the two kinds of grading.    

We follow Proposition \ref{prop:Poisson-like-cohom}   
and apply the same 
 argument as in Example \ref{Eg:Type1} to 
\(\sum_{j=0}^{\dim M} \Lambda^{j}\cbdl{M}\), and get the following fact.   

\begin{kmProp} Assume \(\varphi\) is an even \(p\)-form (i.e., the grade is
\(p+1\)) on a manifold \(M\). 
\begin{enumerate}
\item 
Let \(\varphi\) satisfy \(\Sbt{\varphi}{\varphi} = d (\varphi \we
\varphi) =0 \).  Then \(d_{\varphi}\) satisfies \(d_{\varphi}\circ
d_{\varphi} = 0\) and defines Poisson-like cohomology groups.   
The coboundary operator is \( d_{\varphi} (U) 
= d( \varphi \we U ) \), 
in short, \( d_{\varphi} = d \circ e_{\varphi } = e_{ d \varphi} + e_{ \varphi}
\circ d \) where \( e_{\psi} (U) = \psi \we U\). 

\item 
Let \(\varphi\) be a closed  even \(p\)-form. 
Then \( d_{\varphi} =  e_{ \varphi} \circ d \) and it defines Poisson-like
cohomology groups. 
\kmcomment{
In particular, 
 let \(\varphi\) be a symplectic form. Then
\(d_{\varphi}\) defines Poisson-like cohomology groups. 
}

In particular, 
let \(\varphi\) be a constant function \(1\), which is a 0-form on \(M\). 
The associated Poisson-like cohomology groups
are just the de Rham cohomology groups of \(M\).   
\end{enumerate}
\end{kmProp}

Although 
\( \sum_{p=0}^{\dim M} \Lambda^{p} \cbdl{M} \) for a manifold \(M\) is a
typical superalgebra, its 
chain spaces are infinite dimensional, and it may be hard to manipulate 
those (co)homology groups.   
If \(M\) is a finite dimensional Lie group, then its invariant forms are
finite dimensional, and we would like to handle the superalgebra \(\frakh\) defined by 
\(\frakhN{[-1-j]} = \Lambda^{j} \frakm^{*} \) for \( j=0,\ldots, \dim \frakm \)
where \(\frakm\) is the Lie algebra of a Lie group \(M\) and \(\frakm^{*}\)
is the dual space of \(\frakm\).  The super bracket is given by
\(\Sbt{\alpha}{\beta} = \parity{a} d( \alpha \we \beta )\) where \( \alpha
\in \Lambda^{a} \frakm^{*}\).

If \( \varphi \in \Lambda ^{p} \frakm^{*} \) with \(p\) odd, then 
\(\Sbt{\varphi}{\varphi}=0\) holds automatically but it is not clear that 
\(d _ {\varphi} \circ d _ {\varphi} = 0\) or not, when \( d_{\varphi} (U) =
\Sbt{\varphi}{U} = \parity{p} d( \varphi \we U )\). One concrete example is the
following:
\begin{exam} \label{exam:BadNews}
Start with 4-dimensional Lie algebra of \(\uType[12] 
=\{\Sbt{\yb{1}}{\yb{2}} = 0, 
            \Sbt{\yb{1}}{\yb{3}} = \yb{1}, 
            \Sbt{\yb{1}}{\yb{4}} = -\yb{2}, 
            \Sbt{\yb{2}}{\yb{3}} = \yb{2}, 
            \Sbt{\yb{2}}{\yb{4}} = \yb{1} , 
            \Sbt{\yb{3}}{\yb{4}} = 0 \}\).  
Then \(d\) is given
by  \( 
d( \zb{1}) = -\zb{1}\we \zb{3}-\zb{2}\we \zb{4}\), 
\(
d( \zb{2}) = \zb{1}\we \zb{4}-\zb{2}\we \zb{3}\), 
\( d( \zb{3}) =  d( \zb{4}) = 0\) where \( \zb{i}\) are the dual of \(\yb{j}\)
.   
Take a 1-form \( \varphi = \sum_{i=1}^{4} \cb{i} \zb{i}\). 
The matrix representations of \( d_{\varphi}\) are the followings.   
\[
A(0,2) = 
\left[\begin{array}{cccccc}
0 & \cb{1} & -\cb{2} & \cb{2} & \cb{1} & 0
\end{array}\right], 
\  
A(1,3) = 
\left[\begin{array}{cccc}
2 \cb{2} & 0 & -\cb{4} & \cb{3}
\\
 -2 \cb{1} & 0 & -\cb{3} & -\cb{4}
\\
 0 & 0 & \cb{2} & -\cb{1}
\\
 0 & 0 & \cb{1} & \cb{2}
\end{array}\right], 
\  
A(2,4) = 
\left[\begin{array}{c}
-2 \cb{4} \\
 0 \\
 2 \cb{2} \\
 0 \\
 -2 \cb{1} \\
 0
 \end{array} \right],\] and  \ A(3,5) = A(4,6) = O.  
Since 
\( d_{\varphi}^{2}(\Lambda^{0}\frakm^{*}) =  -2 
\left( \cb{1}^{2}+ \cb{2}^{2}\right) 
\&^{\we} \left(z_{1},z_{2},z_{3},z_{4}\right) \) and 
\( d_{\varphi}^{2}(\Lambda^{j}\frakm^{*})= \) trivial for \(j=1 \ldots 4\).  
Thus, if \( \cb{1}^{2}+ \cb{2}^{2} = 0\) then \( d_{\varphi} \) defines
``cohomology groups''.  
When \( \cb{3} = \cb{4} = 0\) then the coboundary operator is trivial and 
the sequence of Betti number is 
\( [1,4,6,4,1]\), which are \( \binom{4}{j}\) for \( j=0 \ldots 4\). 

If \( \cb{3} \ne 0 \) and \(\cb{4} = 0 \) then non-trivial rank is 
\( \rank A(1,3) = 2\) and the sequence of Betti numbers is [1,2,6,2,1]. 
If \( \cb{4} \ne 0 \)  then non-trivial ranks are  
\( \rank A(1,3) = 2\) and  
\( \rank A(2,4) = 1\),  
so the sequence of Betti numbers is [1,2,5,2,0]. 
\end{exam}

\subsection{Poisson-like cohomology groups by 2-forms of dual Lie algebra}
We study Poisson-like coboundary operator defined by a 2-form \(\varphi\)
satisfying \(\Sbt{\varphi}{\varphi} =0\) on the Lie superalgebra 
\( \sum_{j=0}^{\dim \frakm} \Lambda^{j} \frakm^{*}\)
generated by
the dual space of finite dimensional Lie algebra \(\frakm\).  
Since \( \Sbt{\varphi}{\varphi}
\in \Lambda^{2+2+1}\frakm^{*}\) in general,  if \( \dim \frakm \leqq 4\),
then \(\Sbt{\varphi}{\varphi} =0\) holds automatically. This means every
2-form is a Poisson-like tensor if \( \dim \frakm \leqq 4\). 

\begin{exam} 
Let \(\varphi\) be a 2-form on a dual space of a 4-dimensional Lie algebra. 
Since \(\Sbt{\varphi}{\varphi}\) is 5-form, 
 \(\Sbt{\varphi}{\varphi} = 0\) holds in general.  Thus, every 2-form
 \(\varphi\) defines coboundary operator \(d_{\varphi}\) and Poisson-like
 cohomologies.  
\end{exam}

\newcommand{\cc}{u} 
\begin{exam} Take a 5-dimensional Lie algebra 
\begin{align*}
\uType[8] := \{ & 
 \Sbt{\yb{1}}{ \yb{2}} = 
 \Sbt{\yb{1}}{ \yb{3}} = 
 \Sbt{\yb{1}}{ \yb{4}} = 
 \Sbt{\yb{1}}{ \yb{5}} = 
 \Sbt{\yb{2}}{ \yb{3}} = 
 \Sbt{\yb{2}}{ \yb{4}} = 0 , \\
 & 
 \Sbt{\yb{2}}{ \yb{5}} =  \yb{1} ,
 \Sbt{\yb{3}}{ \yb{4}} = 0  ,
 \Sbt{\yb{3}}{ \yb{5}} =  \yb{3} ,
 \Sbt{\yb{4}}{ \yb{5}} = \cc   \yb{4}\}\; .
\end{align*}
The dual space has the relations 
\[ \left[ d(\zb{1}),\ldots, d(\zb{5})\right] = 
\left[-\zb{2}{\we}\zb{5}, 0, - \zb{3}{\we}\zb{5},
-\cc \zb{4}{\we}\zb{5},0\right]
\]
where \(\cc\) is a real constant.   
Let \(\varphi\) be a 2-form  
\begin{align*}\varphi = &  \cb{1} \zb{1} \we \zb{2}
+ \cb{2 } \zb{1} \we \zb{3}
+ \cb{3 } \zb{1} \we \zb{4}+ \cb{4 } \zb{1} \we \zb{5}
+ \cb{5 } \zb{2} \we \zb{3}
\\ & 
+ \cb{6 } \zb{2} \we \zb{4}
+ \cb{7 } \zb{2} \we \zb{5}
+ \cb{8 } \zb{3} \we \zb{4}
+ \cb{9 } \zb{3} \we \zb{5}
+ \cb{10 } \zb{4} \we \zb{5}
\end{align*} where \( \cb{i}\) are constant for \(i=1,\ldots, 10\). 

\(\Sbt{\varphi}{\varphi} =0\) is equivalent to
\begin{equation}\label{eqn:gen:one}
2(1+\cc)( \cb{1} \cb{8} - \cb{2} \cb{6}  + \cb{3} \cb{5}) = 0\;. 
\end{equation}
If \(1+\cc=0\), any 2-form \(\varphi\) satisfies \(\Sbt{\varphi}{\varphi}=0\)
and defines coboundary operators and Poisson-like cohomologies. 
The matrix representations of \( d_{\varphi}\) are as follows. 
\begin{align*}
A(0,3) &= 
\left[
\begin{array}{*{12}{c}} 
0& 0& 0& 0& \cb{2}& - \cb{3}& 0& \cb{2} + \cb{5}&   \cb{3} - \cb{6}&  0 
\end{array}
\right]
\\
A(1,4) & = 
\left[
\begin{array}{*{5}{c}} 
0& -\cb{5}&  \cb{6}& 0 & -\cb{8} \\ 
0& \cb{2}& - \cb{3}   & 0& 0   \\ 
0& -\cb{1}& 0& 0 & \cb{3}  \\ 
0& 0&  \cb{1} & 0  & - \cb{2}  \\
0& 0& 0& 0& 0 \\
\end{array}
\right]
\\
A(i,i+3) &= 0 \qquad  (i=2,\ldots).  
\end{align*}

If \( 1+\cc \ne 0\), then  \( \cb{1} \cb{8} - \cb{2} \cb{6}  + \cb{3} \cb{5} =
0\) holds if and only if  \(\Sbt{\varphi}{\varphi}=0\). 

If \( \cb{1}= \cb{2} = \cb{3} = \cb{5}= \cb{6} = \cb{8} = 0\) then
\(d_{\varphi} = 0\), i.e.,  \(A(i,i+3) = 0\) for  \(i=0,\ldots\). Otherwise, 
we may assume \( \cb{8} \ne 0\) and so
\( \cb{1} = ( \cb{2} \cb{6}  - \cb{3} \cb{5})/\cb{8}\). 
Then the matrix representations of \(d_{\varphi}\) are as follows.  
\begin{align*}
A(0,3) &= 
\left[\begin{array}{cccccccccc}
0 & 0 & 0 & 0 & \cb{2} & \cc \cb{3}  & 0 & \cb{2}+\cb{5} & \cb{3} + \cc \cb{6} & \cb{8} \left(\cc +1\right) 
\end{array}\right]\;, \quad \rank(A(0,3)) = 1\;, 
\\
A(1,4) &= 
\left[\begin{array}{ccccc}
0 & -\cb{5} & - \cc \cb{6}   & - (\cc+1) \cb{8} & -\cb{8} 
\\
 0 & \cb{2} & \cc \cb{3}   & 0 & -( \cc + 1 ) \cb{8} 
\\
 0 & \ds \frac{-\cb{2} \cb{6}+\cb{3} \cb{5}}{\cb{8}} & 0 & 
  \left(\cc +1\right) \cb{3} & \cb{3}+ (\cc+1)\cb{6} 
\\
 0 & 0 &\ds -\frac{\left(\cb{2} \cb{6}-\cb{3} \cb{5}\right) \cc}{\cb{8}} & 
 -( \cc+1) \cb{2} & -\cb{2} - \left( \cc +1\right) \cb{5}
\\
 0 & 0 & 0 & 0 & 0 
\end{array}\right]\;, \quad \rank(A(1,4)) = 2\;, 
\\
A(2,5) &= 
(\cc+1)
\left[\begin{array}{*{10}{c}}
\cb{8} &  -\cb{6} &  \cb{5} &  0 &  \cb{3} &  -\cb{2} &  0 &  \ds
 \frac{\left(\cb{2} \cb{6}-\cb{3} \cb{5}\right) }{\cb{8}} 
&  0 &  0 
\end{array}\right]^{T}\;, \qquad \rank(A(2,5)) = 1\;. 
\end{align*}

\end{exam}

\section{Dual Poisson-like coboundary operator and cohomologies}
We already handle two kinds of Lie superalgebras, one is 
the direct sum of tangential multi-vector fields with the Schouten bracket and the other is 
the direct sum of differential forms with a Lie super bracket given by \(\ds \Sbt{\alpha}
{\beta} = \parity{a} d( \alpha \we \beta)\) where \( \alpha \) is an \(a-\)form.   
In the context of Poisson-like cohomology theory, we have two long ``exact
sequences'' as follows. 

\adjustbox{scale=0.85, center}{%
\begin{tikzcd}
& \Lambda^{1}\frakm^{} \ar[r, "d_{\pi}"] 
&  \Lambda^{2}\frakm^{} \ar[r,"d_{\pi}"] 
& \cdots  \ar[r,"d_{\pi}"] 
&  \Lambda^{m-1}\frakm^{} \ar[r,"d_{\pi}"] 
&  \Lambda^{m}\frakm^{} \ar[r,"d_{\pi}"] & {0} \\[-5mm]
\Lambda^{0}\frakm^{*} \ar[r, "d"] 
&  \Lambda^{1}\frakm^{*} \ar[r,"d"] 
&  \Lambda^{2}\frakm^{*} \ar[r,"d"] 
& \cdots  \ar[r,"d"] 
&  \Lambda^{m-1}\frakm^{*} \ar[r,"d"] 
&  \Lambda^{m}\frakm^{*} \ar[r,"d"] & {0} 
\end{tikzcd}
}

where \(\ds \frakm = \tbdl{M}\) for a differentiable \(m\)-dimensional
manifold \(M\).   
There may be a duality between multi-vector fields and
differential forms of a finite dimensional Lie algebra      
in some naive sense, we may consider the dual object of
Poisson-like coboundary operators and Poisson-like cohomology groups in the
case.

\begin{exam}
As the first experiment, we handle
the tangential Lie superalgebra associated with 
the 4-dimensional Lie algebra \(\frakm\) which is called Type[1] in
\cite{MR3030954} explained in Example \ref{Eg:Type1},     
and 
the Poisson-like coboundary operator 
which comes from the constant function 1
 explained in Proposition 3.1 (2) on 
the cotangent Lie superalgebra associated with 
the 4-dimensional Lie algebra \(\frakm\).    
Namely, we have two sequences as above. 
In a ``naive'' sense, we may define 
\begin{equation} \label{eqn:Couple}
\langle \Omega , d_{\pi} U \rangle
= \langle  \delta \Omega ,  U \rangle \quad  \text{  for } 
\text{  a \(p\)-form  } \Omega 
\text{ and a \( ( p-1) \)-vector  } 
U 
\end{equation}
where 
\(\langle \cdot, \cdot \rangle\) is the natural pairing. 

We revisit Example \ref{Eg:Type1}.  The Poisson structures are
given by      
\begin{align*}
\pi & = \yb{1} \we ( c_{1} \yb{2}  +c_{2} \yb{3}  +c_{3} \yb{4}  ) 
+c_{4} \yb{2} \we \yb{3} \; ,   
\text{ and } \;   
\pi \we \pi = 
2 c_{3} c_{4} {\&}^{\we} \! \left(y_{1},y_{2},y_{3},y_{4}\right)
\; . \\ 
& d_{\pi} [ \yb{1},\ldots, \yb{4} ] = 
 [d_{\pi} (\yb{1}), \ldots, d_{\pi}(\yb{4})] = 
\left[0,0,-c_{3} y_{1}{\we }y_{2},c_{4} y_{1}{\we}y_{3}
+c_{2} y_{1}{\we}y_{2}\right]
\\
& d_{\pi} [ \yb{1}\we \yb{2}, \ldots, \yb{3} \we \yb{4} ] = 
\left[0,0,0,0,c_{4} {\&^{\we}} \! \left(y_{1},y_{2},y_{3}\right),
-c_{2} {\&^{\we}} \! \left(y_{1},y_{2},y_{3}\right)
-c_{3} {\&^{\we}} \! \left(y_{1},y_{2},y_{4}\right)
\right]
\\
&
d_{\pi} [ \yb{1}\we \yb{2} \we \yb{3}, \ldots, \yb{2} \we \yb{3} \we \yb{4} ] = 
\left[0,0,0,0\right] \; . 
\end{align*}

The Lie algebra structure  defined by 
\( \Sbt{\yb{2}}{\yb{4}} = \yb{1}\), \( \Sbt{\yb{3}}{\yb{4}} = \yb{2}\) induces
the exterior differentiation \(d\) given by 
\( d \zb{1} = - \zb{2} \we \zb{4}\; , \   
 d \zb{2} = - \zb{3} \we \zb{4}\; , \   
 d \zb{3} = 0 \; , \   
 d \zb{4} = 0\; \), where \(\{\; \zb{i} \mid i=1..4\; \} \) are 
the dual basis of  \(\{\; \yb{i} \mid i=1..4\; \} \).   
From \eqref{eqn:Couple}, \(\delta \Omega\) is given by 
\(\ds \delta \Omega = \sum_{U}
\langle \Omega , d_{\pi} U \rangle \text{Dual}\; U \).  
Thus, 
\begin{align*}
& \delta  [ \zb{1}, \ldots,  \zb{4} ] 
= [ \delta(\zb{1}), \ldots, \delta (\zb{4}) ] 
= [0, 0, 0, 0] \\
& \delta  [ \zb{1}\we \zb{2}, \ldots, \zb{3} \we \zb{4} ] = [ -\cb{3}\zb{3} +
\cb{4}\zb{4} , \cb{4}\zb{4}, 0,0,0,0] \\
& \delta  [ \zb{1}\we \zb{2} \we \zb{3}, \ldots, \zb{2} \we \zb{3} \we \zb{4} ] = [ 
\cb{4} \zb{2}\we \zb{4} -
\cb{2}\zb{3} \we \zb{4} , - \cb{3}\zb{3}\we \zb{4}, 0,0] \\
& \delta  [ \zb{1}\we \zb{2} \we  \zb{3} \we \zb{4} ] = [ 0] 
\end{align*}
And we get the following table. 

\[
\begin{array} {c|*{5}{c}} 
 & 0 & 1 & 2 & 3 & 4 \\
 \hline
 \text{
Dim} & 1 & 4 & 6 & 4 & 1 \\
\text{
ImageDim} & 0 & 0 & r &  r & 0 \\
\text{
Betti} & 1 & 4-r & 2(3-r) & 4-r & 1 \end{array}
\]
where \(r = \rank \text{ of 2x2 matrix } \begin{bmatrix} \cb{2} & \cb{3} \\ \cb{4} & 0 \end{bmatrix}
\).  
By taking \(\ds \delta = \delta_{\pi} = \) the dual of \(d_{\pi}\),  we get ``2-dimensional  
cochain complex'' with an anchor the cochain complex of de Rham as below.  

\medskip

\adjustbox{scale=0.85, center}{%
\begin{tikzcd}
   \Lambda^{0}\frakm^{*} \ar[r,"d"]
&  \Lambda^{1}\frakm^{*} \ar[r,"d"]   \ar[ld,  "\delta_{\pi}" ']  
&  \Lambda^{2}\frakm^{*} \ar[r,"d"]   \ar[dl,  "\delta_{\pi}" '] 
&  \Lambda^{3}\frakm^{*} \ar[r,"d"]   \ar[dl,  "\delta_{\pi}" '] 
&  \Lambda^{4}\frakm^{*} \ar[r,"d"]   \ar[ld,  "\delta_{\pi}"']  
& {0} 
\\
   \Lambda^{0}\frakm^{*} \ar[r, "d"] 
&  \Lambda^{1}\frakm^{*} \ar[r,"d"] 
&  \Lambda^{2}\frakm^{*} \ar[r,"d"] 
&  \Lambda^{3}\frakm^{*} \ar[r,"d"] 
&  \Lambda^{4}\frakm^{*} \ar[r,"d"] & {0} 
\end{tikzcd}
}
\bigskip

In this example, \(\delta \circ d=0\) and \( d \circ \delta =0 \) on each
stage.  

\end{exam}

\newcommand{\WE}{\mathit{\&^{\we}}}
\newcommand{\Vol}{\mathit{Vol}}
\newcommand{\kmCC}[2]{{\mathfrak C}^{#1}_{#2}}

\begin{exam}
The 4-dimensional Lie algebra \(\frakm\) Type[2] (
              \( \Sbt{\yb{1}}{\yb{4}} = a \yb{1}, \Sbt{\yb{2}}{\yb{4}} = \yb{2}, 
 \Sbt{\yb{3}}{\yb{4}} = \yb{1}+ \yb{2}\))
in \cite{MR3030954} gives non trivial case. 
The exterior differentiation \(d\) is described as follows. 
\begin{align*}
d([\zb{1},\ldots, \zb{4}]) 
& = [-a \zb{1}  \we  \zb{4}, -\zb{2}  \we  \zb{4} - \zb{3}  \we  \zb{4}, 
    -\zb{3}  \we  \zb{4}, 0]
\\
d([\zb{1}\we\zb{2},\ldots]) 
& = [(a + 1) \WE(\zb{1}, \zb{2}, \zb{4}) + \WE(\zb{1}, \zb{3}, \zb{4}), 
  (a + 1) \WE(\zb{1}, \zb{3}, \zb{4}), 0, 2 \WE(\zb{2}, \zb{3}, \zb{4}), 0, 0]
\\
d([\WE(\zb{1},\zb{2},\zb{3}),\ldots])
& = 
         [(-a - 2) \WE(\zb{1}, \zb{2}, \zb{3}, \zb{4}), 0, 0, 0]
\\
d([\WE(\zb{1},\zb{2},\zb{3}, \zb{4})])
& = 
                              [0]
\end{align*}

A candidate of Poisson tensor is
\(
\pi = \Cb{1} \yb{1}  \we  \yb{2} + \Cb{2} \yb{1}  \we  \yb{3}
       + \Cb{3} \yb{1}  \we  \yb{4} + \Cb{4} \yb{2}  \we  \yb{3}
       + \Cb{5} \yb{2}  \we  \yb{4} + \Cb{6} \yb{3}  \we  \yb{4} 
       \)
for constants 
\(      {\Cb{1}, \Cb{2}, \Cb{3}, \Cb{4}, \Cb{5}, \Cb{6}} \). 
Since
\begin{align*} 
\Sbt{\pi}{\pi} =& 
(-a \Cb{1} \Cb{6} + a \Cb{2} \Cb{5} - \Cb{1} \Cb{6} + \Cb{2} \Cb{5}
   - \Cb{2} \Cb{6} - 2 \Cb{3} \Cb{4}) \WE(\yb{1}, \yb{2}, \yb{3})
   \\& 
   + (a \Cb{3} \Cb{5} - \Cb{3} \Cb{5} - \Cb{3} \Cb{6}) \WE(\yb{1}, \yb{2}, \yb{4})
   \\&
   + (a \Cb{3} \Cb{6} - \Cb{3} \Cb{6}) \WE(\yb{1}, \yb{3}, \yb{4})
   + \Cb{6}^2   \WE(\yb{2}, \yb{3}, \yb{4})
\end{align*}
the Poisson conditions are, 
\begin{align*}
& -a \Cb{1} \Cb{6} + a \Cb{2} \Cb{5} - \Cb{1} \Cb{6}
   + \Cb{2} \Cb{5} - \Cb{2} \Cb{6} - 2 \Cb{3} \Cb{4} = 0, \\ 
& a \Cb{3} \Cb{5} - \Cb{3} \Cb{5} - \Cb{3} \Cb{6} = 0 , \qquad 
  a \Cb{3} \Cb{6} - \Cb{3} \Cb{6} = 0, \qquad  
  \Cb{6}^{2}  = 0 \; .
\end{align*}
Solving the above, we get 5 cases as follows.
\begin{align*}
&\{a = -1, \Cb{1} = \Cb{1}, \Cb{2} = \Cb{2}, \Cb{3} = 0, \Cb{4} = \Cb{4},
\Cb{5} = \Cb{5}, \Cb{6} = 0\}\; , 
\\ & 
 \{a = 1, \Cb{1} = \Cb{1}, \Cb{2} = \Cb{2}, \Cb{3} = \Cb{3}, \Cb{4} =  \Cb{2}
 \Cb{5}/\Cb{3}, \Cb{5} = \Cb{5}, \Cb{6} = 0\}\; , 
 \\& 
  \{a = a, \Cb{1} = \Cb{1}, \Cb{2} = 0, \Cb{3} = 0, \Cb{4} = \Cb{4}, \Cb{5} =
  \Cb{5}, \Cb{6} = 0\}\; , 
  \\ & 
  \{a = a, \Cb{1} = \Cb{1}, \Cb{2} = \Cb{2}, \Cb{3} = 0, \Cb{4} = \Cb{4},
  \Cb{5} = 0, \Cb{6} = 0\}\; , 
  \\ & 
  \{a = a, \Cb{1} = \Cb{1}, \Cb{2} = \Cb{2}, \Cb{3} = \Cb{3}, \Cb{4} = 0,
  \Cb{5} = 0, \Cb{6} = 0\}\; . 
  \end{align*}

Case 1: \( \pi =  \Cb{1} \yb{1}  \we  \yb{2} + \Cb{2} \yb{1}  \we  \yb{3}
       + \Cb{4} \yb{2}  \we  \yb{3} + \Cb{5} \yb{2}  \we  \yb{4}\), 
       where \(
      {a = -1, \Cb{3} = 0, \Cb{6} = 0}\). 
\(\pi\we \pi =  -2 \Cb{2} \Cb{5} \WE(\yb{1}, \yb{2},
  \yb{3}, \yb{4}) \). 
Poisson coboundary operator \(d_{\pi}\) is 
\begin{align*}
d_{\pi}[ \yb{1},\ldots, \yb{4}] 
& = [-\Cb{5} \yb{1}  \we  \yb{2}, 0, -\Cb{5} \yb{2}  \we  \yb{3}, \Cb{2} \yb{1}  \we  \yb{2} + 2 \Cb{4} \yb{2}  \we  \yb{3} + \Cb{5} \yb{2}  \we  \yb{4}]
\\
d_{\pi}[ \yb{1}\we \yb{2}, \ldots] 
& = 
[0, 0, -2 \Cb{5} \WE(\yb{1}, \yb{2}, \yb{4}) - 2 \Cb{4} \WE(\yb{1}, \yb{2}, \yb{3}), 0, 0, -\Cb{2} \WE(\yb{1}, \yb{2}, \yb{3})]
\\
d_{\pi}[\WE( \yb{1},\yb{2}, \yb{3}),\ldots ] 
& = [0, 0, -\Cb{5} \WE(\yb{1}, \yb{2}, \yb{3}, \yb{4}), 0] \; . 
\end{align*}
The dual operator \(\delta\) of \(d_{\pi}\) is given by 
\begin{align*}
\delta[ \WE(\zb{1},\ldots,\zb{4})]  
&  =   [-\Cb{5} \WE(\zb{1}, \zb{3}, \zb{4})]
\\
\delta[ \WE(\zb{1},\zb{2},\zb{3}),\ldots]  
& = 
      [-2 \Cb{4} \zb{1}  \we  \zb{4} - \zb{3}  \we  \zb{4} \Cb{2}, 
        -2 \Cb{5} \zb{1}  \we  \zb{4}, 0, 0]
\\
\delta[ \zb{1}\we \zb{2},\ldots]  
&= 
[\Cb{2} \zb{4} - \Cb{5} \zb{1}, 0, 0, 2 \Cb{4} \zb{4} - \Cb{5} \zb{3}, \Cb{5} \zb{4}, 0]
\\
\delta[ \zb{1},\ldots,\zb{4}]   
& = [0, 0, 0, 0]
\end{align*}
We verify \(\delta \circ \delta = 0\) directly from the above relations.  
In this case, \(d \circ \delta\) and \( \delta \circ d\) are not trivial and satisfy  
\begin{alignat*}{5}
(d \circ \delta) \text{ and } (\delta \circ d) \text{ of } [ \WE( \zb{1},\ldots, \zb{4}) ] 
& =                             [0] \text{ and } [0] \\
(d \circ \delta) \text{ and } (\delta \circ d)  \text{ of }  [ \WE( \zb{1},\zb{2}, \zb{3}), \ldots ] 
& = [0, 0, 0, 0] \text{ and }
              [\Cb{5} \WE(\zb{1}, \zb{3}, \zb{4}), 0, 0, 0]
\\
(d \circ \delta) \text{ and } (\delta \circ d)  \text{ of }   [ \zb{1}\we \zb{2}, \ldots ] 
& = 
  [-\Cb{5} \zb{1}  \we  \zb{4}, 0, 0, \Cb{5} \zb{3}  \we  \zb{4}, 0, 0]
\text{ and }
 [0, 0, 0, 0, 0, 0]
\\
(d \circ \delta) \text{ and } (\delta \circ d)  \text{ of }   [ \zb{1},\zb{2} , \ldots ] 
& = [0, 0, 0, 0] \text{ and } [0, -\Cb{5} \zb{4}, 0, 0] \; .
\end{alignat*}

We leave the handling of other cases to readers.  Case 2, 3 and 5 are non-trivial,
namely  \( \delta \circ d \ne 0\) and \( d \circ \delta \ne  0\) in some stage.  
Even \(\delta\) in Case 4 is non-trivial but 
\( \delta \circ d =0\) and \( d \circ \delta = 0\).  
\begin{remark}
We constructed the dual \(\delta\) of Poisson cohomology operator
\(d_{\pi}\) formally. 
It may be  interesting to study the geometric meanings of 
\( \delta \circ d \) or \( d \circ \delta \) for the original Poisson
structure \(\pi\). 
\end{remark}

\end{exam}

\subsection{Multi-vector fields and differential forms with polynomial coefficients on \(\mR^{n}\) }

As explained several times, two
prototypes of Lie superalgebras are 
\begin{enumerate}
\item \(\sum_{i=1}^{\dim M} \Lambda^{i} \tbdl{M} \) with the Schouten bracket
\(\SbtS{\cdot}{\cdot}\) 
(cf.\ \cite{Mik:Miz:super2}), and  
\item \(\sum_{i=0}^{\dim M} \Lambda^{i} \cbdl{M} \) with 
\( \Sbt{ \alpha}{\beta} = \parity{a} d( \alpha \wedge \beta ) \) where \( \alpha
\in \Lambda^{a}\cbdl{M}\) (cf.\ \cite{Mik:Miz:superForms}). 
\end{enumerate}

As concrete examples of  finite dimensional sub superalgebras, we know  
the space of 
multivector fields with homogeneous polynomials as coefficient on 
$n$-dimensional number space \(\mR^{n}\) with the 
Cartesian coordinates \(\xb{1},\ldots, \xb{n}\) by 
the Schouten bracket (cf.\ \cite{Mik:Miz:super2}), 
and the differential forms with homogeneous polynomials with the modified
super bracket by the exterior differentiation \(d\).  
We denote \( d( \xb{i} ) = \zb{i}\), then \( d( \zb{i}) =0\) holds and  
\(d (P) = \sum_{i=1}^{n} \frac{\pdel P}{\pdel \xb{i}} \zb{i} \) for a
polynomial \(P\), and we see the following.  
\begin{equation} 
\label{eqn:d} 
d ( \kmCC{m}{k} ) \subset \kmCC{m+1}{k-1}
\quad \text{where}\quad  
\kmCC{m}{k} = m\text{-forms with } k\text{-homogeneous polynomials as
coefficient} \; .
\end{equation}
In the same way, we set 
\(\ds \kmC{k}{m} = m\text{-vector fields with } k\text{-homogeneous polynomials}\). 
\begin{kmProp}
\begin{equation}
\dim \kmCC{m}{k}\ \text{ and }\   
\dim \kmC{k}{m}\quad \text {are both equal to }\  \tbinom{n}{m} \tbinom{n-1+k}{n-1}. 
\end{equation}
Thus, 
\( \kmC{k}{m} = \{0\}\) and \(\kmCC{m}{k} = \{0\}\) when \( m > n \), and 
\(\kmC{k}{m} = \{0\}\) and \(\kmCC{m}{k} = \{0\}\) when \( k < 0 \).  
\end{kmProp}
\begin{remark}
The relation \eqref{eqn:d} shows  
the sum \( m+k\) of two indices \(m\) and \(k\) of \(\kmCC{m}{k}\) is
``invariant'' under the operation \(d\). 
\end{remark}
\begin{kmProp}
Let \(\pi\in \kmC{h}{2}\) with  \(\Sbt{\pi}{\pi} = 0 \), namely  
\(h\)-homogeneous Poisson 2-vector field on \(\mR^{n}\), and let 
\(d_{\pi} = \SbtS{\pi}{\cdot}\).  It holds that 
\begin{equation} \label{eqn:d:pi}
d_{\pi} ( \kmC{k}{m} ) \subset \kmC{h + k-1}{m+1}\; .
\end{equation}
\end{kmProp}
\begin{remark}
Looking for some ``invariant'' by two indices \(m,k\) of \(\kmC{k}{m}\)  
in \eqref{eqn:d:pi}, 
the sum \( m+k\) is invariant when \(h=0\) and 
the difference \(m-k\) is invariant when \(h=2\).  
\end{remark}
\begin{definition}
Take a 2-homogeneous Poisson 2-tensor field \(\pi\) on \(\mR^{n}\). 
Let \( \Vol = \WE( \zb{1}, \ldots, \zb{n})\) (the standard volume form). 
By the natural pairing \(   \langle \Vol, U \rangle \) of \(U \in \kmC{k}{m}\)
and \(\Vol\), we get \(\langle \Vol, U\rangle \in \kmCC{n-m}{k}\) and get the
dual \(\delta\) of \(d_{\pi}\) as below. 

\medskip 

\adjustbox{scale=0.95, center}{%
\begin{tikzcd}
\kmC{k}{m}\ni U \ar[r, "d_{\pi}"] \arrow[d, "vol"'] & \SbtS{\pi}{U}\in \kmC{k+1}{m+1}
\ar[d, "vol"] \\
\kmCC{n-m}{k}\ni \langle \Vol,U\rangle \arrow[u] \ar[r, dashed, "\delta"'] & \langle \Vol, \SbtS{\pi}{U}
\rangle\in \kmCC{n-(m+1)}{k+1} \arrow[u]
\end{tikzcd}
}

Thus, there is  a chain complex: 
\begin{equation} 
\kmCC{n}{k} \mathop{\longrightarrow}^{\delta}  
\kmCC{n-1}{k+1} \mathop{\longrightarrow}^{\delta} \ldots \mathop{\longrightarrow}^{\delta}  \kmCC{1}{k+n-1}   \mathop{\longrightarrow}^{\delta}  
\kmCC{0}{k+n} \;. 
\end{equation}

\end{definition}

Since \( \kmC{2}{2}\) is spanned by the product of 
\(\{ \xb{i} \xb{j} \mid  1 \leqq i\leqq j \leqq n\}\) and 
\(\{ \yb{i} \we \yb{j} \mid  1 \leqq i < j \leqq n\}\),  
a general 2-homogeneous 2-vector field \(\pi\) is given by 
\begin{equation} 
\pi = \sum_{i\leqq j,k < \ell} C[i,j:k,\ell ] \xb{i} \xb{j} \yb{k} \we \yb{\ell} 
\quad \text{
with}\quad\tbinom{n}{2}\tbinom{n-1+2}{2}\;\text{-parameters}\quad C[i,j:k,\ell ]\ .  \end{equation}
The Poisson conditions are given by 
\begin{align*}
\Sbt{\pi}{\pi} =& 
\sum \sum  C[i,j:k,\ell ] 
 C[i',j':k',\ell{}'] \Sbt{ 
\xb{i} \xb{j} \yb{k} \we \yb{\ell}} {\xb{i'} \xb{j'} \yb{k'} \we \yb{\ell '}}
\\ =& 
\sum \sum  C[i,j:k,\ell ] 
C[i',j':k',\ell{}'] \left( 
\xb{i}\xb{j}\Sbt{\yb{\ell}}{\xb{i'}\xb{j'}} \WE (\yb{k},\yb{k'},\yb{\ell '}) 
-\xb{i}\xb{j}\Sbt{\yb{k}}{\xb{i'}\xb{j'}} \WE (\yb{\ell},\yb{k'},\yb{\ell '}) 
\right. 
\\ 
& \left. \hspace{55mm}
-\xb{i'}\xb{j'}\Sbt{\yb{k'}}{\xb{i}\xb{j}} \WE (\yb{k},\yb{\ell},\yb{\ell '}) 
+\xb{i'}\xb{j'}\Sbt{\yb{\ell'}}{\xb{i}\xb{j}} \WE (\yb{k},\yb{\ell},\yb{k '}) 
\right)\;. 
\end{align*}
For instance, let \(n=4\). Then we need 60 parameters \(C[i,j:k,\ell]\). Now
we denote them in 1-dimensional way like \(C[1],\ldots, C[60]\). The 
Poisson conditions are given by 80 bi-homogeneous equations. 
Additionally, we care about the rank of Poisson structure. 
The rank is 4 if \( \pi \we \pi \ne 0\) in this case, and is 2 if 
 \( \pi \we \pi = 0\) and \( \pi  \ne 0\). 
Still we have gotten the equations, 
it is hard to show the all of them  here,  nor to solve the given equations by hand, so we ask help by Symbol calculus. 
It provides general solutions of 51 types.  One of them is the following.

The independent variables are 
\(
\{\CSp{4}, \CSp{7}, \CSp{10}, \CSp{14}, \CSp{37}, \CSp{40}\}
\). 
The dependencies are 
\begin{align*}
& \CSp{1} = 0, \CSp{2} = 0, \CSp{3} = 0, \CSp{5} = 0, \CSp{6} = 0, \CSp{8} = 0
, \CSp{9} = 0, \CSp{11} = 0, \CSp{12} = 0, \CSp{13} = 0, \CSp{15} = 0, \CSp{16} = 0
,  \\
& \CSp{17} = \frac{\CSp{7} \left(\CSp{14}-\CSp{37}\right)}{\CSp{4}}, \CSp{18} = 0
, \CSp{19} = 0, \CSp{20} = 
-\frac{\CSp{7} \CSp{40}-\CSp{10} \CSp{14}+\CSp{10} \CSp{37}}{\CSp{4}},
\CSp{21} = 0, 
\CSp{22} = 0, \CSp{23} = 0,\\
& \CSp{24} = 0, \CSp{25} = 0, \CSp{26} = 0, \CSp{27} = 0
, \CSp{28} = 0, \CSp{29} = 0, \CSp{30} = 0, \CSp{31} = 0, \CSp{32} = 0, \CSp{33}
 = 0, \CSp{34} = 0, \CSp{35} = 0,\\
 & \CSp{36} = 0,
  \CSp{38} = 0, \CSp{39} = 0, 
\CSp{41} = 0, \CSp{42} = 0, \CSp{43} = 0, \CSp{44} = 0, \CSp{45} = 0, \CSp{46} = 0
, \CSp{47} = 0, \CSp{48} = 0, \CSp{49} = 0, \\ 
& \CSp{50} = 0, \CSp{51} = 0, \CSp{52} = 0, \CSp{53} = 0, \CSp{54} = 0, \CSp{55} = 0, \CSp{56} = 0, \CSp{57} = 0, 
\CSp{58} = 0, \CSp{59} = 0, \CSp{60} = -\CSp{37}\; .  
\end{align*}

Assume \( \CSp{4}=\CSp{37}=1, \CSp{10}=\CSp{14}=\CSp{40}=0\). 
Then the Poisson tensor is 
\begin{equation}\label{eqn:myTgt}  \pi = 
+     (  \xb{1}
+ \CSp{7}\xb{2})\xb{4} \yb{1} \we \yb{2}
- \CSp{7}\xb{2}\xb{4}  \yb{1} \we \yb{3}
+        \xb{2}\xb{4}  \yb{2} \we \yb{3}
-        \xb{4}^2      \yb{3} \we \yb{4} \; . 
\end{equation}
\begin{remark}
When given \eqref{eqn:myTgt},  \(\Sbt{\pi}{\pi}=0\) 
and \( \pi \we \pi = -2 ( \xb{1} + \CSp{7} \xb{2}) \xb{4}^3  
\WE ( \yb{1}, \yb{2}, \yb{3}, \yb{4}) \) 
are straightforward  reconfirmations.  
\end{remark}

We only show the table of the Poisson cohomology groups derived 
from the Poisson structure \eqref{eqn:myTgt} the following left. 
\[
\begin{array}{c|*{9}{c}} 
    & \kmC{0}{1}   & \rightarrow 
    & \kmC{1}{2}   & \rightarrow 
    & \kmC{2}{3}   & \rightarrow 
    & \kmC{3}{4} \\\hline
\text{Dim} & 4 && 24 && 40 && 20 \\\hline
\text{Rank} & 3 && 15 && 10 && 0 \\\hline
\text{Betti} & 1 && 6 && 15 && 10 
\end{array}
\qquad 
\qquad 
\begin{array}{c|*{10}{c}} 
    & \kmCC{3}{0}   & \rightarrow 
    & \kmCC{2}{1}   & \rightarrow 
    & \kmCC{1}{2}   & \rightarrow 
    & \kmCC{0}{3} 
    \\\hline
\text{Dim} &   4 && 24 && 40 && 20 \\\hline
\text{Rank}&   3 && 15 && 10 && 0 \\\hline
\text{Betti}&  1 && 6 && 15 && 10 
\end{array}
\]
The table right above shows the behavior of the dual of the Poisson cohomology
with respect the volume form \(\Vol = \WE( \zb{1},..,\zb{4})\).  
There is the de Rham cochain complex by \(d\) on 
\(\kmCC{0}{3}\rightarrow \kmCC{1}{2}\rightarrow \kmCC{2}{1}\rightarrow \kmCC{3}{0}\), 
which is always trivial in the sense that the each Betti number is 0.  And
there is a ``double'' complex as the following.  It may be interesting to
study how the
double complex shows properties of 2-homogeneous Poisson 2-vector field. 

\bigskip

\adjustbox{scale=1.10,center}{%
\begin{tikzcd}
   \kmCC{0}{3} \ar[r,"d"]
&  \kmCC{1}{2} \ar[r,"d"]   \ar[ld,  "\delta_{\pi}" ']  
&  \kmCC{2}{1} \ar[r,"d"]   \ar[dl,  "\delta_{\pi}" '] 
&  \kmCC{3}{0} \ar[r,"d"]   \ar[dl,  "\delta_{\pi}" '] 
& {0} 
\\
   \kmCC{0}{3} \ar[r, "d"] 
&  \kmCC{1}{2} \ar[r,"d"] 
&  \kmCC{2}{1} \ar[r,"d"] 
&  \kmCC{3}{0} \ar[r,"d"] 
& {0}
\end{tikzcd}
}

\nocite{Mik:Miz:super3}
\nocite{Mik:Miz:Sato:Engel}

\printbibliography

\end{document}
\end{document} \endinput

\nocite{Mik:Miz:super2}
\nocite{Mik:Miz:super3}
\nocite{Mik:Miz:Sato:Engel}

\printbibliography

\end{document}

%% file: km_cmd.tex
\newcommand{\kmqed}{\hfill\ensuremath{\blacksquare}}

\newcommand{\rank}{\operatorname{rank}}

\newcommand{\ds}{\ensuremath{\displaystyle }}
\newcommand{\pdel}{\partial}

\newcommand{\mR}{\ensuremath{\mathbb{R}}} 
\newcommand{\mZ}{\ensuremath{\mathbb{Z}}} 
\newcommand{\mH}{\ensuremath{\mathbb{H}}} 

\newcommand{\ANYb}[2]{{#1}_{#2}}

\newcommand{\xb}[1]{x_{#1}}
\newcommand{\yb}[1]{y_{#1}}
\newcommand{\zb}[1]{z_{#1}}

\newcommand{\frakg}{\mathfrak{g}}
\newcommand{\frakgN}[1]{{\mathfrak{g}}_{#1}}

\newcommand{\frakh}{\mathfrak{h}}
\newcommand{\frakhN}[1]{\mathfrak{h}_{#1}}

\newcommand{\tbdl}[1]{\mathrm{T}(#1)}

\newcommand{\cbdl}[1]{\mathrm{T}^{*}(#1)}

\newcommand{\Sbt}[2]{[#1,#2]}                
\newcommand{\SbtS}[2]{[#1,#2]_{\text{\footnotesize S}}} 
                
\newcommand{\parity}[1]{(-1)^{#1}}

\newcommand{\Bkt}[2]{\{#1,#2\}_{\pi}}

\newcommand{\kmC}[2]{ C_{#1}^{#2}}       

\newcommand{\CSp}[1]{\text{C}_{#1}} 

\newcommand{\we}{\wedge}



%% file: km_cmd2.tex
\renewcommand{\dim}{\textrm{dim}}

\renewcommand{\[}{$$} \renewcommand{\]}{$$}



%% file: km-env-e.tex
\usepackage{ntheorem} 
{%
\theoremheaderfont{\bfseries}
\theorembodyfont{\rmfamily}

\newtheorem{defn}{\textbf{Definition}}

\newtheorem{prop}{Proposition}[section]
\newtheorem{exam}{Example}[section]
\newtheorem{remark}{Remark}[section]

\newtheorem{definition}{Definition}[section]
\newtheorem{theorem}{Theorem}[section] 
\newtheorem{thm}{Theorem}[section]
\newtheorem*{thm-none}{Theorem}[section]
\newtheorem{kmProp}[theorem]{Proposition}

\newtheorem{kmCor}[theorem]{Corollary}

}

\renewcommand{\[}{$$} \renewcommand{\]}{$$}

\newcommand{\frakm}{\mathfrak{m}}